\documentclass[12pt]{article}

\usepackage{amssymb, latexsym, amsmath}
\newtheorem{Theorem}{Theorem}[section]
\newtheorem{Proposition}[Theorem]{Proposition}

\newtheorem{Corollary}[Theorem]{Corollary}
\newtheorem{Remark}[Theorem]{Remark}
\newtheorem{Definition}[Theorem]{Definition}
\newtheorem{Question}[Theorem]{Question}
\newtheorem{Conjecture}[Theorem]{Conjecture}
\newcommand{\bi}{\begin{enumerate}}
\newcommand{\ei}{\end{enumerate}}
\newcommand{\be}{\begin{equation}}
\newcommand{\ee}{\end{equation}}
\newcommand{\ba}{\begin{array}}
\newcommand{\ea}{\end{array}}

\def\c{\mathbb{C}}
\def\cp{\mathbb{CP}}

\def\r{\mathbb{R}}

\def\z{\mathbb{Z}}

\def\T{\mathbf{T}}
\def\t{\mathbf{t}}
\def\th{\theta}

\author{ Massimiliano Pontecorvo}
\title{Bi-Hermitian and locally conformally K\"ahler surfaces %\\ Kato surfaces and bi-Hermitian structures 
}
\begin{document}
\maketitle %\today

%\begin{abstract}
%\end{abstract}

%\input intro

\section{Introduction}

We intend to report on some work which has been going on for about twenty years on bi-Hermitian surfaces, by which we mean Riemannian four--manifolds $(M^4,g)$ equipped with two \it integrable \rm complex structures $J_+$ and $J_-$ which are both $g$-orthogonal inducing the same orientation. 

\begin{Definition} An (almost) complex structure $J\in End(TM)$ with $J^2=-id.$ on an even-dimensional Riemannian manifold $(M,g)$ is said to be orthogonal if 
$g(JX,JY)=g(X,Y)$ for all tangent vectors $X,Y\in TM$. 
\end{Definition}

It is important to stress here that in the framework of bi-Hermitian metrics each $J_\pm$ is assumed to be integrable, in the usual sense that induces local holomorphic coordinates. This condition imposes  constrains on the curvature of the Riemannian metric $g$; in fact, by conformal invariance, it imposes conditions on the Weil tensor $W$ of the metric.

The following \it local \rm result of Salamon provided the initial motivation for studying bi-Hermitian metrics in four-dimension:

\begin{Theorem}\cite{sa94}
Each point $p\in (M^4,g,or.)$ of an oriented Riemannian four-manifold belongs to a neighborhood in which there are either zero, one, two or infinitely many complex structures (up to sign) which are $g-$orthogonal, inducing the given orientation. %The last alternative occurs if and only if $g$ is anti-self-dual.
\end{Theorem}

Notice that the last alternative occurs if and only if $g$ is anti-self-dual. Soon after the above local result the following question appeared to be quite natural:

\begin{Question}  \bf Salamon 90's. \it What about Riemannian four-manifolds admitting exactly two (distinct, up to sign) orthogonal complex structures inducing the same orientation ? 
\hspace{0.5 cm} \rm In what follows, such a structure $(g,J_\pm)$ will be called a bi-Hermitian metric. 
\end{Question}

The first local results are due to Piotr Kobak \cite{ko99} who explicitly constructed bi-Hermitian metrics on compact tori.
Then came the remarkable work of Apostolov-Gauduchon-Grantcharov \cite{agg99} where the local theory was developed along with global constructions coming from deformations of hyperk\"ahler metrics. % on tori and K3 surfaces.

\noindent The present paper is devoted to discuss the following \bf global version \rm of Salamon question: 
\begin{Question} 

\noindent a) Which compact four-manifolds support bi-Hermitian metrics ?  \\
b) More precisely, which compact complex surfaces $(S,J_+)$ admit a Riemannian metric $g$ and a \rm different \it complex structure $J_-$ such that $J_\pm$ are $g$-orthogonal and  induce the same orientation ?
\end{Question}

\begin{Remark} \rm
The precise meaning of \it different \rm in this question is that there is $p\in S$ where, as endomorphisms of the tangent bundle we have, $J_+(p) \neq \pm J_-(p)$. It will then turn up that the two complex structures may or may not be biholomorphic, in any case they will be in the same class of Enriques-Kodaira classification. This is the framework we intend to work in, because it soon appeared that in order to answer question b) is enough to assume $(S,J_+)$ to be a \it minimal \rm complex surface \cite{fp10} \cite{cg11}. This means, by Castelnuovo criterium, free from smooth rational curves of self-intersection $-1$.
\hfill  $\bigtriangleup$  \end{Remark}

We end up this introduction with a quick look at dimensions greater than four.  %\medskip

\noindent \bf High dimensions. \rm
We mention here briefly about two possible higher dimensional generalizations of this setting in the Riemannian and in the quaternionic  case.
The Weil tensor $W$ of a \it Riemann \rm metric becomes irreducible as soon as the  real dimension $n\geq 5$ and certainly the condition $W=0$ implies, for $n$ is even, the existence of an abundance of Orthogonal Complex Structures, simply because the twistor space is a complex manifold. Conversely, Salamon showed \cite{sa95} that each orthogonal complex structure imposes linear conditions on the conformal curvature $W$, which however has a big number of components, difficult to control.

A different behavior was observed in the higher dimensional \it quaternionic \rm setting where the natural analogue of an oriented conformal structure $(M^4,[g],or.)$ is a $GL(n,\mathbb H)Sp(1)$-structure on the tangent bundle of $M^{4n}$. It was shown in \cite{amp99} that for every quaternionic dimension $4k\geq 8$ the local existence of two compatible complex structures forces the almost quaternionic structure to be a quaternionic structure.
By the integrability of the corresponding twistor space, this is equivalent to say that in this situation
we always have infinitely many complex structures compatible with the quaternionic structure, locally. 
We conclude that there is no analogue notion of bi-Hermitian metric in this setting.

\smallskip \noindent
\it Acknowlegments. \rm Most of our results have been obtained in collaboration with Prof. Akira Fujiki which we hearth-fully thank for the longstanding collaboration and friendship. It is a great pleasure to say that I have been learning a lot about the subject from Simon Salamon, Piotr Kobak, Gueo Grantcharov, Paul Gauduchon, Vestislav Apostolov, Georges Dloussky, Ryushi Goto and Liviu Ornea.

\section{Fundamental divisor and flat bundles}

\subsection{The fundamental line bundle of a bi-Hermitian metric.}

In what follows we will denote by $S$ the compact complex surface $S=(M,J_+)$ where $M$ is a four-manifold with bi-Hermitian structure given by two complex structures $J_\pm$ which are orthogonal to the same Riemannian metric and induce the same orientation. 

Our goal is to give informations about the complex geometry of $S$.
%and althought it is not possible to give a complete classification of all all such $S$ in the spirit of the Enriques-Kodaira classification, their smooth structure is certainly well understood. 
To this end we start to introduce here what will be called the \it fundamental divisor \rm of $S$. It is the effective divisor $T$ geometrically described below as the intersection of smooth complex surfaces in the almost complex, real $6$-dimensional twistor space $Z$ of $(M,[g],or.)$. The space $Z$ only depends on the conformal class of the metric which we denote by $[g]$.

The standard notation $[D]$ for the holomorphic line bundle of the divisor $D$ will be used here together with additive notation for divisors and multiplicative for line bundle so that $[D+D']=[D]\otimes [D']$. 
We start by saying that the \it support \rm $\mathcal T$ of the fundamental divisor $T$ is easily seen, by the description below, to be the set of points where the two complex structures agree, up to sign:
$$\mathcal T = \{p\in M \;|\; J_+(p)=\pm J_-(p) \}.$$

As a matter of fact $\mathcal T$ turns out to be the closure of a union of smooth 2-dimensional real submanifolds  \cite{fp14}. However, being a non-necessarily transverse intersection, the divisor $T$ will typically be singular and highly non-reduced, see \cite{fp15} for a computation of all its multiplicities.

Setting $\mathbf{t}:=b_0(\mathcal T)$ it follows from surface classification that the number of connected components can only be 
\begin{equation}\label{components}  \mathbf{t}=0,1,2 
\end{equation}
in other words $T$ can have at most two connected components and may well be empty.

We end up this subsection with some considerations involving the twistor space.
% by just saying that a Chern class computation in the twistor space $Z$ of the oriented four manifold % yields the following: 

The bi-Hermitian structure $([g],J_\pm)$ defines two complex surfaces in the twistor space $Z=\{J\in \operatorname{End}(TM) \; | \; J^2=-id.\}$ of all (locally defined) almost complex structures compatible with the metric and orientation of $M$, for which we will use the following notation $S_\pm:=J_\pm(M)\subset Z$. Similarly we set $X_\pm:=S_\pm\amalg \sigma(S_\pm)$ where $\sigma:Z\to Z$ is the involution $J\mapsto -J$. Notice that there is a tautological biholomorphism $(M,J_\pm)\cong S_\pm$.

In this notation $T\subset S=(M,J_+)$ can be identified with either of the following two subsets in the twistor space: $X_+\cap S_-$ or  $S_+\cap X_-\subset Z$; this exhibits $T$ as an almost complex subvariety of either $S_+$ or $S_-$; $T$ is therefore a complex curve in each of the two smooth surfaces, in particular closed in the analytic Zariski topology and nowhere dense. %by \cite[1.3]{po97}.
A first Chern class computation in $Z$ then yields the following: 

\begin{Proposition}\label{feq} Let $t:Z\to M$ denote the twistor projection, the image  $t(T)=t(X_+\cap S_-)=t(S_+\cap X_+) \subset M$ has a natural structure
of divisor $T$ in both surfaces $S_\pm$ with the property
 $$c_1(T)= c_1(S_\pm)$$
For each compact complex surface $S_\pm$ there is a holomorphic line bundle $F_\pm$ such that
\begin{enumerate}
\item \quad $c_1(F_\pm)=0$ 
\item  \quad $[T] = K^{-1}_\pm \otimes F_\pm$ where $K_\pm$ is the canonical line bundle of the surface $S_\pm$.
\end{enumerate}
\end{Proposition}

\begin{Definition} We will call $T$ the \rm fundamental \it divisor of the bi-Hermitian surface and $F=[T]\otimes K$ the \rm fundamental flat line bundle. \it
\end{Definition}

We summarize this discussion with the following statement which gives an important obstruction to the existence of bi-Hermitian surfaces

\begin{Corollary}
For every bi-Hermitian surface $S=(M,J_+)$ the fundamental divisor $T$ is what is called a numerically anti-canonical (NAC) divisor:  simply because its Chern class is opposite to that of the canonical bundle: $[T] =F\otimes K_S^{-1}$ with $c_1(F)=0$.
\end{Corollary}

\subsection{ Lee bundle of a LCK metric}

For dimensional reasons, the fundamental $(1,1)$-form $\omega$ of an Hermitian surface \newline $S=(M^4,g,J)$ always satisfies the equation
$$d\omega =\omega \wedge \th$$
for a unique $1$-form $\th$ called the \it Lee form \rm of the Hermitian metric such that:
%. Then, the following holds
\bi
\item The metric is K\"ahler if and only if $\th=0$ 
%\item the Lee form $\theta ^\prime$ of a conformally related metric $g^\prime=e^{-f}\,g \;$ is given by $\theta^\prime = \t -df$
%, we easily see that:
%\bi
 \item The Lee form $\th=df$ is exact when the metric $e^{-f}\,g$ is K\"ahler -- i.e. $(g,J)$ is a globally conformally K\"ahler metric.
 \item The Lee form is closed -- i.e. locally exact -- precisely when the Hermitian metric is conformal to K\"ahler, locally. We abbreviate this condition, which was introduced by Vaisman in the 70's and is still a topic of active research in every complex dimension \cite{ov24b}, by LCK metric.
 \ei

\begin{Definition} We say that a compact complex Hermitian surface is LCK if it admits a LCK metric which is not globally conformally K\"ahler.
\end{Definition}

In the above situation, the closed Lee form $\th$ defines a flat line bundle $L_{\th}\in \check{H}^1(M,\r^{>0})$ which we call here the \it Lee bundle \rm -- or also called  \it weight bundle \rm by Gauduchon et al. -- with local trivializations $e^{-(f_i)}$ and transition functions given by $\{ g_{ij}=\frac{e^{-(f_i)}}{e^{-(f_j)}} \}$ which are \it positive constants \rm because the locally exact Lee form satisfies $df_i=\th=df_j$ on double intersections while the cocycle condition on triple intersections holds automatically.

\begin{Remark} \rm
These same considerations show that the correspondence between closed Lee forms and Lee bundles extends to a linear map in cohomology 
$$\exp : H^1_{dR}(S) \to \check{H}^1(S,\r^{>0})$$ 
from the first deRham cohomology group to the multiplicative group of line bundles with positive constant transition functions (and therefore flat) on $M$.  %% -- under tensor product $\otimes$ -- 
\hfill  $\bigtriangleup$       \end{Remark}

For the purpose of the present work, the most important instance of this correspondence is given by the following result of  Appostolov-Gauduchon-Grantcharov:

\begin{Proposition}\cite{agg99}[Lemma 1, Prop. 2] For any bi-Hermitian metric  $(M^4,g,J_\pm)$ consider the respective two Lee forms $\theta_\pm$ of $(g,J_\pm)$ which of course are real and need not be closed individually. However, the sum $\theta_+ + \theta_-$ is always closed and corresponds to the fundamental flat bundle $F$ of the bi-Hermitian structure: $F:= L_{\frac12 (\theta_+ + \theta_-)}$.
\end{Proposition}

%\bigskip\bigskip  \noindent  %\it 
\subsection{Gauduchon degree}

\begin{Remark} % The Gauduchon metric and the Gauduchon degree . Degree of flat line bundles \rm \label{gauduchon metric}
\rm By a strong result of Gauduchon \cite{ga84},  \it in every conformal class $c$ of Hermitian metrics on a compact complex surface $S$ there is a unique representative $g\in c$ whose Lee form $\theta$ is coclosed or equivalently whose fundamental $(1,1)$-form $\omega$ is $\partial\overline\partial$-closed.
Furthermore the Lee form is coexact -- $\th = \delta G$ -- if and only if the second Betti number is even \cite{ga84}. \rm

Such a Riemannian metric $g$ will be called the \it Gauduchon metric \rm of the conformal class and this result has many useful applications among which we mention the notion of \it degree \rm of a holomorphic line bundle on the Hermitian surface $S=(M,g,J)$. 

Because we will only be interested in \it the special case of flat line bundles \rm $L_a \in H^1(M,\r^{>0})$ associated to closed $1$-form $a\in H^1_{dR}(M)$ the formula in the definition simplifies to \cite{ad16} 
 \be 2\pi \operatorname{deg}(\mathcal L_a)=-\langle a,\theta \rangle = 
        -\langle a_h,\theta_h \rangle 
     \label{scalarproduct}  \ee
where $\theta$ is the Lee form of the Gauduchon metric $g$ with $g$-harmonic part denoted by $\theta_h$ while $a_h$ denotes the $g$-harmonic part of the $1$-form $a$ and $\langle , \rangle$ is the global inner product with respect to $g$; the last equality holds because $a$ is closed and $\theta$ is coclosed.

We will make the following use of degree: let $\mathcal L$ denote the holomorphic line bundle which is the complexification of $L$, then $\operatorname{deg}$ measures the $g$-volume with sign of a virtual meromorphic section of $\mathcal L$ -- even if $\mathcal L$ might not admit non-trivial meromorphic sections.  In particular $deg (\mathcal L)>0$ whenever  $H^0(S,\mathcal L)\neq 0$ and $\mathcal L \neq \mathcal O$.
Finally, because the space of Gauduchon metrics is connected and $\r\setminus {0}$ is not, the sign of $deg(L)$ is independent of the particular choice of Gauduchon metric on $(S,J)$. It may, however, depend on $J$ \cite{fp15}.
\hfill  $\bigtriangleup$  \end{Remark}

The above mentioned result of Apostolov-Grantcharov-Gauduchon allows to compute the degree of the fundamental flat line bundle with respect to the $J_+$--Gauduchon metric, see \cite{ap01} %metric is 
$$2\pi \operatorname{deg}_{J_{+}} (F) = -\langle \theta_+ , \frac12 (\theta_+ +\theta_-) \rangle = -\frac14 ||\theta_+ +\theta_- ||^2$$  because the equality $||\theta_+||^2 = ||\theta_- ||^2$ holds for any metric in the conformal class \cite{po97}. 
From this we get the result that the fundamental flat bundle $F$ of a bi-Hermitian surface has
\begin{equation}\label{degree}
\operatorname{deg}(F)\leq 0 
\end{equation}
with equality if and only if $F$ is trivial

We end this section by observing that the Gauduchon degree of a Lee bundle $\mathcal L _\theta$ --- i.e. the holomorphic flat line bundle of a LCK metric --- is always strictly negative, because from equation (\ref{scalarproduct}) it follows that $\operatorname{deg} (\mathcal L _{\theta}) = -||\theta_h||^2 < 0$ otherwise $g$ would be conformally K\"ahler and $b_1(S)$ even. 
% In particular then, $H^0(S,\mathcal L _{\theta})=0$ and the same holds for any positive power.

%\input bi-Hermitian
\section{bi-Hermitian surfaces}

Getting back to the main question: which compact complex surfaces $(S,J)$ admit a bi-Hermitian structure $(g,J_\pm)$ with $J=J_+$? The aim of this section is to present all surfaces which are known to admit such metrics. It is convenient to spilt our presentation according to the parity of the first Betti number which, in complex dimension $2$, is equivalent to the existence of a K\"ahler metric.

\subsection{K\"aherian case.}

It turns out to be convenient to start with the case $b_1(S)$ even, because of the following
\begin{Proposition}\cite[Lemma 4]{agg99} When $(M^4,g,J_\pm)$ is compact bi-Hermitian with \rm even \it
first Betti number the two Lee forms always satisfy the condition that 
\be   \textrm{the sum } \quad  \theta_+ + \theta_- \quad\text{ is exact. } \ee      \label{gke}
\end{Proposition}
\it Proof. \rm Let $g$ be the $J_+$-Gauduchon metric, as $b_1(M)$ is even, we can assume $\theta_+ =\delta\alpha$. On the other hand $d(\theta_+ + \theta_-) =0$ so that the orthogonal decomposition of $\theta_-$ into harmonic, coexact and coclosed part is of the form $\theta_-=\theta^h_- - \delta\alpha +df$. Because the two Lee forms have the same global $L^2$-norm we conclude that $\theta^h_-=0=df$ so that $\theta_++\theta_-=0$ and $g$ is $J_-$-Gauduchon too.
 \hfill$\Box$\medskip

The fundamental equation (\ref{feq}) reduces to
\be \T=K^{-1}  \label{f=0}   \ee
and the necessary condition to admit bi-Hermitian metrics is then $H^0(S,K^{-1})\neq 0$.
 
As a consequence, the canonical bundle of $S$ is trivial or else $S$ must be a ruled surface. In the first case, explicit examples where constructed by Kobak \cite{ko99} and Apostolov-Grantcharov-Gauduchon \cite{agg99} on tori and K3 surfaces, while Enriques or hyperelliptic surfaces cannot be bi-Hermitian because the flat bundle $F$ is not real.
 
The remaining case in also of very high interest because the above proposition says that bi-Hermitian surfaces are indeed \it generalized K\"ahler \rm for $b_1$ even.
Important contributions by Hitchin \cite{hi06}, Gualtieri \cite{gu14}, Goto \cite{go12} et al. proved the converse statement that indeed any ruled surface with anti-canonical divisor is bi-Hermitian, when $b_1$ is even.

\subsection{Non-K\"ahlerian surfaces.}

This case is where we gave a contribution and where some questions are still open. When $b_1(S)$ is odd, the flat line bundle $F$ may not be trivial, equivalently the closed 1-form $\theta_+ + \theta_-$ may or may not be exact. In any case $S$ cannot be a surface of Kodaira class VI:

\begin{Proposition}\cite{ap01} \label{negdeg}
The canonical bundle of a compact bi-Hermitian surface $S$ of odd first Betti number has negative degree, as a consequences all plurigenera must vanish and $Kod(S)=-\infty$. In other words the surface $S$ belongs to class VII of Kodaira classification.% and in particular $b_1(S)=1$.
\end{Proposition}
\it Proof. \rm
By the fundamental equation $T=K^{-1}\otimes F$ and recall that $T$ is effective or zero while $F$ has negative degree or is trivial \ref{negdeg}. It follows that deg $K = $ deg $F $- vol $T <0$ unelse $F$ and $T$ are both trivial. However, since $b_1(S)$ is odd, a result of Apostolov-Gauduchon-Grantcharov \cite[Prop. 4]{agg99} says that when $F$ is trivial the fundamental divisor $T$ must have two non-empty connected components, therefore cannot be trivial.
\hfill$\Box$\medskip

From now on the bi-Hermitian surface $S$ belongs to class VII$_0$ which exactly means: Kod$(S)=-\infty$, $S$ is minimal with odd first Betti number. Recall that indeed all such surfaces satisfy $b_1(S)=1$ and have negative definite intersection form so that the topological signature $\sigma(S)=-b_2(S)$ agrees with the Chern number $c_1^2(S)=2\chi(S) +3\sigma(S)=-b_2(S)$.

Although there is no complete classification of surfaces in this class, we now present a few results which are particularly useful for discussing bi-Hermitian surfaces. We start from the case of vanishing first Chern class.
 
\subsubsection {Surfaces in class VII$_0$ and $b_2(S)=0$}
A Theorem of Bogomolov \cite{bg76}, later clarified by Teleman \cite{te94} and Yau et al. \cite{lyz94} asserts that there exist only two different types of surfaces $S\in$ VII$_0$ with $b_2(S)=0$. 

The first class is given by Hopf surfaces which are defined to be compact quotients of $\c^2\setminus\{0\}$; they have fundamental group isomorphic to $\z\oplus(\z/\z_{k})$ for $k\geq 1$ and are diffeomorphic to $S^1\times S^3$, up to finite coverings.. The easiest example is when the action is generated by the contraction $(z,w)\mapsto \frac12 (z,w)$. The complete list of contractions was described by Kodaira and shows that 
 anti-canonical divisors are always effective, consisting of two elliptic curves in the diagonal case and of one curve with multiplicity in the other case. 

It follows that every Hopf surface satisfies the necessary condition (\ref{feq}) for existence of bi-Hermitian metrics.  %%%

Indeed, it has been shown in \cite{po97} and \cite{ad08} that Hopf surfaces admit an abundance of bi-Hermitian structures, with all possible values of $\t = b_0(T)=0,1,2$. 

\medskip
By Bogomolov theorem the only other case satisfying $S\in$VII and $b_2(S)=0$ is that $S$ is a so called Inoue-Bombieri surface \cite{in74} \cite{bm73}. Such surfaces are always discrete quotients of $\mathcal H \times\c$, the product of the upper-half plane with a complex line; they have no curves at all.
Inoue-Bombieri surfaces however cannot admit bi-Hermitian metrics by \ref{negdeg} because the degree of their canonical bundle is always positive by a result of Teleman \cite{te06}.
%\end{document}
\subsubsection{Surfaces in class VII$_0^+$}

In order to investigate bi-Hermitian surfaces we are now left with the non-K\"ahlerian case of negative Kodaira dimension and positive second Betti number: Kod$(S)=-\infty$ and $b_2(S)>0$. Recall furthermore that we can assume that $S$ is \it minimal; \rm the standard notation for this class of surfaces is VII$_0^+$. 
%meaning that $S$ is not obtained by blowing-up points on another surface; equivalently, by Castelnuovo's criterion %$S$ admits no smooth rational curve of self-intersection $-1$.
Althought the Kodaira classification of these surfaces is still an open problem, thanks to the work of Dloussky all the bi-Hermitian surfaces in this class belong to a well known class of surfaces which will be presented in the last section. 

\begin{Proposition} \cite{dl06} Every bi-Hermitian surface $S\in $ VII$_0^+$ is a Kato surface.
\end{Proposition}
\it Proof. \rm The fundamental divisor $T$ of any bi-Hermitian surface is either effective or zero with the same Chern class of $S$: $c_1(T)=c_1(K^{-1}_S)=c_1(S)$ so that $c_1^2(S)=b_2(S)>0$ implies that $T$ is a non-trivial NAC (numerical (pluri-) anti-canonical) divisor on $S\in$ VII$_0^+$ which must therefore be a Kato surface by a strong result of Dloussky \cite{dl06}.
\hfill$\Box$\medskip

\section{Bi-hermitian and LCK metrics, an interplay}

The first examples of bi-Hermitian Kato surfaces, are due to Lebrun \cite{le91} on so called Parabolic Kato surfaces. Their construction is a byproduct of the LeBrun ansatz for self-dual metrics on the connected sum of arbitrary many copies of $\cp_2$. As such these bi-Hermitian metrics are also LCK, by a result of \cite{bo86}.

This circle of ideas motivated the construction of Fujiki-Pontecorvo which was aiming at producing LCK metrics on Kato surfaces. Although the relevant log-deformation theory of singular twistor spaces turned out to be obstructed in the general case, it was proved in \cite{fp10} that the construction produces bi-Hermitian surfaces on all Kato surfaces admitting a disconnected anti-canonical divisor, of real type. This result completed the classification of bi-Hermitian Kato surfaces under the hypothesis that the flat line bundle $F$ is actually trivial.

\begin{Theorem}\cite{fp10}
Every Hyperbolic Inoue surface and every Parabolic Inoue surface of real type admits bi-Hermitan metrics and these are the only possibilities in class VII$_0^+$, when $\t=2$.
\end{Theorem}

Indeed, coming from twistor theory, these metrics are automatically anti-self-dual and therefore LCK. Indeed a strong motivation for this work was to give new inputs to the following old interesting 

\begin{Question}\bf Vaisman 1970's \rm
Which compact complex surfaces admit LCK metrics?  % see \cite{po14} for an overview.
\end{Question}

We refer to  \cite{po14} for an overview of results and just mention here that a consequence of \cite{fp10} was the first proof that all Kato surfaces with Dloussky number equal to $2b$ or $3b$ are LCK; see the next section for the definition of Dloussky number. 

Soon after, this result was greatly generalized by the remarkable work of Marco Brunella which was able to show that \it all \rm Kato surfaces are LCK \cite{br11}.

So that, after Brunella's work, the state of the art about Vaisman question is: all \it known \rm compact complex surfaces are LCK, except for the Inoue-Bombieri surfaces of Belgun \cite{be00}.

%However, a different very interesting connection between bi-Hemitian and LCK geometry of surfaces was later on produced by the work of Apopstolov-Bailey-%Dloussky, 

A different, powerful interplay between LCK and bi-Hermitian metrics is given by the following result of Apostolov-Bailey-Dloussky which, in the same way as LCK geometry can be thought of as a twisted version of K\"ahler geometry, appears as a twisted version of a Hamiltonian construction for generalized K\"ahler metrics of Gualtieri \cite{gu10}.

\begin{Theorem}\cite{abd15} Let $S=(M^4,J)$ be a compact complex surface in class VII with a locally conformally K\"ahler metric and denote by $F$ its Lee bundle.  Then, $M$ also admits a bi-Hermitian metric $(g,J_\pm)$ whose fundamental divisor is $T$ and $J=J_+$ if and only if there is a \rm connected \it NAC divisor $T$ of index 1 -- i.e.  $[T]=K^{-1}\otimes F$ -- with the same real flat bundle $F$ of the LCK metric
\end{Theorem}

This theorem was used by Apostolov-Bailey-Dloussky to produce the first bi-Hermitian metrics on Hopf surfaces with $\t=1$ \cite{abd15}; and then by Fujiki-Pontecorvo to prove the same result for parabolic Inoue surfaces \cite{fp19}.  In fact, surprisingly enough, the Lee bundle of the LCK bi-Hermitian metric turns out to be exactly what is needed for the Apostolov-Bailey-Dloussky theorem to work. The conclusion is that Hopf surfaces admit bi-Hermitian metrics with all possible values of $\t=0,1,2$ while Parabolic Inoue are bi-Hermitian with $\t=1,2$; it is not hard to see that these are the only surfaces which can possibly be bi-Hermitian for different values of $\t$.

\section{Kato surfaces}

The classification of compact complex surfaces in class VII$_0^+$ has received much attention in the last forty years from the Japanese school of Nakamura, Inoue, Enoki as well from the French-based school of Dloussky, Oeljeklaus, Toma and Teleman; we present some of their results in this section. 
For the moment however, the only known surfaces in this class can be described by the following geometric construction of Kato \cite{ka77}, who also showed that they can be deformed into blown-up Hopf surfaces and are therefore always diffeomorphic to $(S^1\times S^3)\# b \overline{\cp}_2$.

A \it Kato surface \rm $S$ is geometrically obtained by successive blow ups of the unit ball in $\c^2$ starting from the origin and continuing to blow up a point on the last constructed exceptional curve until, in order to get a \it minimal \rm compact complex surface, one takes the quotient by a biholomorphism $\sigma$ which identifies a small ball centered at a point $p$ in the last exceptional curve with a small ball around the origin $0\in\c^2$ in such a way that $\sigma(p)=0$, see the picture on \cite[p.5]{dl14}. This process constructs $b:=b_2(S)$ rational curves $D_1,\dots , D_b$ of self-intersection $\leq -2$ which are the only rational curves in $S$.
 
 In other words $S=S_{\pi,\sigma}$ is obtained by identifying, via $\sigma\circ\pi$, the two boundary components of $\hat B \setminus \sigma(B)$ --- denoted by $\partial \hat B \; \amalg \; \sigma(\partial B)$ --- each of which is a $3$-sphere. 
This shows that $S$ has a global spherical shell, briefly called GSS, and $\pi_1(S)=\z$. 

% \bigskip
\begin{Conjecture}
The \bf GSS conjecture \rm of Nakamura \cite{na89} asserts that every surface in class VII$_0^+$ is a Kato surface.
\end{Conjecture}

\medskip
As a matter of fact it is not even known if every surface in class VII$_+^0$ must have curves and, since the algebraic dimension vanishes, these surfaces can have at most finitely many curves.
 However, by a very relevant result of Dloussky-Oeljekluas-Toma \cite{dot03} a surface $S\in$ VII$_0^+$ with $b_2(S)$ rational curves must be a Kato surface and the configuration of the rational curves is encoded by the Dloussky sequence Dl$(S)$ which we now present.

The blow-up construction of  $S$ produces exactly $b_2(S)=:b$ rational curves most of which will have self-intersection number $-2$ because a curve $C$ with self-intersection number $C^2=-(k+2)\leq -3$ can only be obtained by repeatedly blow-up at a fixed node of previously created exceptional curves and will therefore come along with a chain of $-2$-curves of length $(k-1)$.
All other curves in $S$ are either obtained in this way or else by blowing up a general point of the previously created exceptional divisor, in particular their
self-intersection number is $-2$. Therefore, if $(a_1,\dots ,a_b)$ denotes the string of  \it opposite \rm self-intersection numbers of all rational curves in $S$ we can separate it into \it singular \rm sequences $s_k:=((k+2)2\dots 2)$ of length $k$ and \it regular \rm sequences $r_m:=(2\dots 2)$ of length $m$.

This notation was introduced by Dloussky in \cite{dl84} and the \it Dloussky sequence \rm of a Kato surface $S$ is a string of $b=b_2(S)$ positive integers describing the opposite self-intersection number of the rational curves, following their \it order of creation, \rm and grouped into singular sequences $s_k$ and regular sequences of maximal length $r_m$.
\bigskip

We refer to \cite{fp15} for a picture of the dual graph of curves in a Kato surface $S$ with 
%% We consider in this section a Kato surface $S$ associated to a 
\it simple \rm Dloussky sequence
%% as described in \cite[p.335]{ot08} or \cite[p.43]{dl11} of the following form
\be   \mathrm{DlS}=[s_{k_1}s_{k_2}...s_{k_{p}}r_m]  \label{simple} \ee
%with $k_1 \geq 1$ and $m \geq 1$.

A weaker biholomorphic invariant of Kato surfaces is the Dloussky number $\sigma(S)$ which is defined to be the sum of opposite self-intersection numbers of the rational curves $D_1,\dots , D_b$ in $S$ so that $\sigma(S):=-\sum_{i=1}^n D_i^2$. It follows from the previous discussion that 
\begin{equation} 2b \leq \sigma(S) \leq 3b
\end{equation}

We now give a quick presentation of these surfaces using the terminology of the Japanese school.
Kato surfaces with $2b = \sigma(S)$ -- or equivalently Dl$(S)=[r_b]$ -- are called Enoki surfaces and are geometrically characterized by the property of having a cycle $C$, of $b$ rational curves and self-intersection $C^2=0$, which therefore represents a flat line bundle. The versal family of Enoki surfaces is $b$-dimensional and has a $1$-dimensional subfamily of surfaces containing an elliptic curve. These special Enoki surfaces are called Parabolic Inoue. By the previously mentioned results, every Parabolic Inoue surface with line bundle $[C]$ of real type has bi-Hermitian metrics with $\t=2$ as well as different bi-Hermitian metrics with $\t=1$. To the contrary, a general Enoki surface cannot be bi-Hermitian simply because there are no anti-canonical divisors.

The other extreme case $\sigma(S) = 3b$ is when all blow-ups occur at intersection points so that Dl$(S)=[s_{k_1}s_{k_2}...s_{k_{p}}]$ presents no regular sequences. As a consequence, every rational curve belongs to a cycle and there are at most two cycles in $S$. The corresponding surface is then called Hyperbolic Inoue when there are two cycles ($p$ is even). Every such surface is bi-Hermitian by \cite{fp10}, while in the other case $p=$odd there is no bi-Hermitian metric because the only NAC divisor is not of real type.

\bigskip
All other Kato surfaces satisfy the strict inequalities $2b < \sigma(S) < 3b$ and are therefore called \it intermediate. \rm In this case the union of all rational curves is connected and contains a unique cycle with branches departing from it. The number of branches equals the number of regular sequences in Dl$(S)$.

There are no examples of bi-Hermitian intermediate Kato surfaces but in fact there are some serious known obstructions which we now would like to present.

The first obstruction is the index. This is defined to be the least positive integer for which $H^0(S,K^{-m}\otimes F)\neq 0$, for some flat line bundle $F\in$ Pic$_0(S)$.

\smallskip
Therefore, if an intermediate Kato surface $S$ is bi-Hermitian it needs to have index 1.
%\smallskip

\noindent This (index=1-)condition only depends on Dl$(S)$ and as such poses severe restrictions on the configuration of curves, like for example that the number of irreducible components in the cycle must be bigger than that of the branches, see \cite{fp15}.

More obstructions come from the complex structure of intermediate Kato surfaces. Before getting into the relevant complex deformations theory, we need to mention the following convenient way of Dloussky-Oeljeklaus to parametrize flat line bundles. %% on Kato surfaces.  

%================

%\medskip
%The GSS property also gives the following topological result that will be used later.
Following \cite{do99}[\S 1.2], for any Kato surface $S$ the group Pic$_0(S)$ of flat holomorphic line bundles is non-compact and isomorphic to $\c^*$.
Furthermore, being a surface with a global spherical shell, $S$ has a Leray cover by two open sets whose intersection has two connected components (like in the circle) so that every $F\in$ Pic$_0(S)$ is determined by a pair of transition functions $(1,e^{f})$ with $f$ a complex constant. The notation used to indicate this correspondence is as follows

$$\begin{array}{ccl}
    \operatorname{Pic}_0(S)   & 	\longleftrightarrow & \c^* \\
           F  & \longmapsto     &   L^{e^{f}}
 \end{array}$$
and notice that  \it real \rm flat line bundles have $f\in\r$ while $f=0$ if and only if $F$ is trivial.

\medskip
Now we come to \it logarithmic deformations. \rm  More precisely, we consider the versal family $\mathcal S$ of deformations of the pair $(S,D)$ where $S$ is intermediate Kato surface and $D$ is the union of all its rational curves. $D$ is always a connected and reduced normal crossing divisor with $b$ irreducible components. The complex deformation family of $(S,D)$ is unobstructed of complex dimension $h^1(S,\Theta(\log (-D)))=m$, where $m$ is the total length of the regular sequences in Dl$(S)$. Notice that the deformation preserves the divisor $D$ and even more so Dl$(S)$. In particular, there exists an integer $k\geq 2$ which can be defined in various ways \cite{dl14} \cite{fp15} and explicitly computed as a polynomial of degree $p$ in the variables $(k_1,\dots , k_p)$. So that $k$ remains fixed in any logarithmic deformation.  

The following result says that every intermediate Kato surface $S$ has a unique numerically anti-canonical divisor. Furthermore, $S$ is deformation equivalent to a surface $S_u$ with a holomorphic vector field and Dl$(S_u)=$Dl$(S)$, if and only if the index is 1. Finally, recall that the sign of degree is independent of the Gauduchon metric.

\begin{Theorem}\cite{do99}[4.5 and 5.2]\label{exp}
Let $\mathcal S \to U$ be a versal family of deformations of the pair $(S,D)$ with (fixed) integer $k\geq 2$, as above.
%of an intermediate Kato surface with fixed Dloussky, let $k=k($Dl$(S))\geq 2$ denote the exponent of the second component in the Favre germ or equivalently the twisting exponent..
Then, there is a surjective holomorphic map $\mu:U\to\c^*$ such that $\mu(u)$ is the unique complex number for which $H^0(S_u,K^{-1}\otimes L^{\mu(u)})\neq 0$. Furthermore, $\mu$ can be chosen to satisfy the following properties:
\begin{enumerate}
\item The flat bundle $L^{\mu(u)}$ is of real type if and only if $\mu(u) \in \r^{>0}$; 
  the sign of its degree is positive for $\mu(u)>1$ and therefore negative when $0 < \mu(u) < 1$.
\item $S_u$ has a holomorphic anti-canonical section if and only if index$(S)=1=\mu(u)$.
\item $S_u$ has a holomorphic vector field if and only if index$(S)=1$ and $\mu(u) =\frac 1k$.
%\item deg$ L^{\mu(u)})>0$ for $\mu(u)>1$ and deg$L^{\mu(u)})<0$ for $0 < \mu(u) < 1$.
\end{enumerate}
\end{Theorem}

The above versal family can be described using contracting germ of the biholomorphism $\sigma$ in the GSS construction of the Kato surface $S$. This germ is a polynomial map $\phi:\c^2 \to \c^2$ which has been put in the following normal form by Favre \cite{fa00}.
\begin{equation}\label{germ}
\phi(z,\xi)=(\lambda \xi^s z + Q(\xi) , \xi^k)
\end{equation}

The degree of the second component is exactly the integer $k\geq 2$ of the Theorem and can be shown to regulate existence of vector fields. 
The function $\mu$ depends on the first coefficient $\lambda\in\c^*$ of the germ. 

The relevance of the germ technique and the role of the integer $k$ shows up strongly in the following recent result of Apostolov-Dloussky who used Favre germ (\ref{germ}) to construct a negative plurisubharmonic function on the universal cover of any intermediate Kato surface (of arbitrary index) with automorphic factor exactly equal to $k$. In particular, by their general result \cite{ad23}[3.7], $L^{\mu(u)}$ cannot be the Lee bundle of a LCK metric on $S_u$, unless $\mu(u) > k$.

This is an obstruction for the Lee bundle of a LCK metric on $S$ %% but recall that any Kato surface is LCK by Brunella.
which induces an obstruction for bi-Hermitian metrics via the Apostolov-Bailey-Dloussky Theorem \cite{abd15}.

%Finally, as noticed by Apostolov-Dloussky \cite{ad23}[Theorem 4.3] we also get an 
%obstruction for bi-Hermitian metrics because the Apostolov-Bailey-Theorem \cite{abd13} gives a necessary and %sufficient condition when the topological invariant $\t=1$, always true for intermediate surfaces.
 \smallskip
We close the paper with some simple concluding remarks.
Every compact bi-Hermitian four-manifolds with $b_1(M)$ is always diffeomorphic to a blown up Hopf surface $(S^1 \times S^3) \# b \overline{\cp}_2$. However the finer question of determining the complex structure remains open in the intermediate case, only: 

\begin{Question} Can intermediate Kato surfaces support bi-Hermitian metrics ?
\end{Question}

As we have seen, this question is very much related to the study of possible Lee classes of LCK metrics on intermediate Kato surfaces.
This is in fact a very interesting question, and an active field of research for LCK manifolds of any dimension \cite{ov24a}. In complex dimension 2 it has been extended by Apostolov-Dloussky to the study of locally conformally symplectic structures, see \cite{ad16}, \cite{ad23}.

\smallskip
We conclude with a final remark about Vaisman question: first of all, the answer is not always positive, as pointed out by Belgun \cite{be00}. Indeed, in his thesis at \`Ecole Polytechnique he showed that some of the Inoue-Bombieri surfaces are not LCK, proving at the same time that the LCK property is not stable under deformations of the complex structure.

A simple remark is then that a solution of the GSS conjecture would also imply that Belgun examples are the only possible non-LCK surfaces and this would provide a complete answer to Vaisman question. This has been observed for example by Ornea-Verbitsky-Vuletescu in \cite{ovv21}.

\newcommand{\bysame}{\leavevmode\hbox to3em{\hrulefill}\,}

%\end{document}

\end{document}